\documentclass[conference]{IEEEtran}
\IEEEoverridecommandlockouts

\usepackage[utf8]{inputenc}
\usepackage{amsmath}  
\usepackage{amssymb}   
\usepackage{amsthm}
\usepackage{amsfonts}
\usepackage{bbm}
\usepackage{mathtools}
\usepackage{geometry}
\newtheorem{Theorem}{Theorem}
\newtheorem{Lemma}{Lemma}
\newtheorem{Remark}{Remark}
\newtheorem*{Proof}{Proof}

\usepackage{xcolor}
\usepackage[normalem]{ulem}
\usepackage{hyperref}
\usepackage{url}
\usepackage[maxbibnames=99, backend=biber]{biblatex}
\usepackage{amsmath,amssymb,amsfonts}
\usepackage{algorithmic}
\usepackage{graphicx}
\usepackage{textcomp}
\usepackage{xcolor}
\def\BibTeX{{\rm B\kern-.05em{\sc i\kern-.025em b}\kern-.08em
    T\kern-.1667em\lower.7ex\hbox{E}\kern-.125emX}}

\begin{document}

\title{Asymptotics for power posterior mean estimation}

\author{
\IEEEauthorblockN{Ruchira Ray}
\IEEEauthorblockA{\textit{Department of Statistics} \\
\textit{Columbia University}\\
New York, United States \\
ruchira.ray@columbia.edu}
\and
\IEEEauthorblockN{Marco Avella Medina}
\IEEEauthorblockA{\textit{Department of Statistics} \\
\textit{Columbia University}\\
New York, United States \\
marco.avella@columbia.edu }
\and
\IEEEauthorblockN{Cynthia Rush}
\IEEEauthorblockA{\textit{Department of Statistics} \\
\textit{Columbia University}\\
New York, United States \\
cynthia.rush@columbia.edu}
}
\maketitle

\begin{abstract}
Power posteriors “robustify” standard Bayesian inference by raising the likelihood to a constant fractional power, effectively downweighting its influence in the calculation of the posterior. Power posteriors have been shown to be more robust to model misspecification than standard posteriors in many settings. Previous work has shown that power posteriors derived from low-dimensional, parametric locally asymptotically normal models are asymptotically normal (Bernstein-von Mises) even under model misspecification. We extend these results to show that the power posterior moments converge to those of the limiting normal distribution suggested by the Bernstein-von Mises theorem. We then use this result to show that the mean of the power posterior, a point estimator, is asymptotically equivalent to the maximum likelihood estimator.
\end{abstract}

\begin{IEEEkeywords}
power posterior, $\alpha$-posterior, Bernstein-von Mises, posterior mean convergence
\end{IEEEkeywords}

\section{Introduction}
In this work, we study $\alpha$-posteriors (also known as power posteriors), which are proportional to the product of the prior and the likelihood raised to a constant, fractional power \cite[Chapter 8.6]{GhoshalVdV}. This procedure (also known as posterior tempering) decreases the influence of the likelihood in the calculation of the posterior. Decreasing the influence of the likelihood has been explored for several problems as a way to improve Bayesian procedures \cite{M03, WH01, Z06}. More recent literature about this area -- under the name $\alpha$-posteriors -- has focused on the robustness to model misspecification and data corruption of the $\alpha$-posterior compared to that of the ordinary posterior both empirically \cite{SafeBayes, GvO2017, MD19} and theoretically \cite{MD19, MORV22}. To this end, both the theoretical and empirical properties of this procedure have been investigated in recent years \cite{SafeBayes, MD19, BPY19, AR20, MORV22}. In particular, \cite{MORV22} established that $\alpha$-posteriors derived from low-dimensional, parametric models are asymptotically normal. In other words, they obey a Bernstein-von Mises theorem. For standard posteriors, the Bernstein-von Mises theorem says that the large sample distribution of the posterior is close to a normal distribution centered around the maximum likelihood estimator (MLE) with posterior variance depending on the ``curvature” of the likelihood function \cite{TPE, KV12, Miller21}. The work in \cite{MORV22} shows that raising the likelihood to the power $\alpha$ while calculating the posterior distribution leaves the mean of the limiting Gaussian unchanged but divides the variance of the limiting Gaussian by $\alpha$.

We extend the results of \cite{MORV22} by showing that the $k^{th}$ moment of the $\alpha$-posterior converges in probability to the $k^{th}$ moment of the limiting normal distribution suggested by the Bernstein-von Mises theorem in \cite{MORV22}. This result (Theorem \ref{thm:moments_alt}) holds under similar conditions to those of \cite{MORV22}, but we additionally assume the finiteness of the $k^{th}$ $\alpha$-posterior moment to prove $k^{th}$ moment convergence. We use Theorem \ref{thm:moments_alt} to establish the $\sqrt{n}$-asymptotic normality of the $\alpha$-posterior mean (Theorem \ref{thm:BayesEstimator}). In particular, we show that the $\alpha$-posterior mean has the same asymptotic distribution as that of the MLE. This result shows that, asymptotically, tempering does not affect point estimation from tempered posteriors.

The proof of our first result is adapted from the Bernstein-von Mises proofs of \cite{KV12} and \cite{MORV22} in order to show convergence of the moments of the $\alpha$-posterior rather than convergence to a limiting normal distribution in total variation. The second result follows almost immediately from our first result. The main assumptions required in our proofs are slightly stronger than those from \cite{MORV22} and \cite{KV12}, as discussed in \eqref{sec:assumptions}. Like \cite{MORV22} and \cite{KV12}, our results make no assumptions of model correctness or independent and identically distributed (i.i.d.) data.

\section{Main Results}
\subsection{Notation and preliminaries}
Let $\phi(\cdot|\mu,\Sigma)$ denote the multivariate normal density with mean $\mu$ and covariance matrix $\Sigma$. The indicator function of an event $A$ is denoted $\mathbbm{1}\{A\}$. For a sequence of distributions $P_n$ on random variables (or matrices) $X_n$, if $\lim_n P_n(\lvert|X_n\rvert|_2 >\epsilon ) = 0$ for every $\epsilon > 0$, we say $X_n=o_{P_n}(1)$. We say $X_n$ converges in $L_1$ to a random variable $X$ if $\lim_{n\rightarrow\infty}\mathbb{E}[|X_n-X|]=0$ and $\mathbb{E}[|X_n|]$ and $\mathbb{E}[|X|]$ exist. We say  $X_n$ is `bounded in $P_n$-probability', or $O_{P_n}(1)$, if for every $\epsilon > 0$ there exists $M_\epsilon > 0$ such that $P_n(\lvert|X_n\rvert|_2<M_\epsilon) \geq 1-\epsilon$. We use $||X||_p = \left(\sum_i|X_i|^p\right)^{1/p}$ to denote the $p$-norm. Let $B_{v}(\delta)$ denote a closed ball of radius $\delta$ around vector $v \in \mathbb{R}^p$. 

\textbf{Statistical Model:} Let $\mathcal{F}_n = \{f_n(\cdot|\theta): \theta \in \Theta \subseteq \mathbb{R}^p \}$ be a parametric family of densities used as a statistical model for the random sample $X^n = (X_1,\dots,X_n)$. This model is allowed to be misspecified  in the sense that $f_{0,n}$, the true density of the random sample $X^n$, may not belong to $\mathcal{F}_n$.  We will use the notion of the pseudo-true parameter $\theta^*$, which is the true data-generating $\theta$ if the model is well-specified, meaning $f_{0,n}(X^n|\theta)\in\mathcal{F}_n$. If the model is misspecified, $\theta^*$ is the value of $\theta\in\Theta$ such that the model $f_n(X^n|\theta^*)$ is closest in KL-divergence to the true model. Denote the (pseudo) MLE by $\hat{\theta}^{ML} = \arg\max_\theta f_n(X^n|\theta)$, where the asymptotic behavior of the MLE is determined with respect to the model $f_{0,n}$. Let $\mathbb{E}_{f_{0,n}}$ and $\mathbb{P}_{f_{0,n}}$ denote the expectation and probability taken with respect to $f_{0,n}$.

\textbf{Power posterior:} Given the statistical model $\mathcal{F}_n$, a prior density $\pi(\theta)$ for $\theta$, and a scalar $\alpha > 0$, the power posterior, or $\alpha$-posterior, has the density:
\begin{equation}\label{alphaposterior}
    \pi_{n,\alpha}(\theta|X^n) = \frac{f_n(X^n|\theta)^\alpha\pi(\theta)}{\int_\Theta f_n(X^n|\theta)^\alpha\pi(\theta) d\theta}.
\end{equation}

\subsection{Assumptions}\label{sec:assumptions}
Our main results hold under the following assumptions:
\begin{itemize}
    \item[\textbf{(A0)}] The MLE, $\hat{\theta}^{ML}$, is unique and there exists a $\theta^*$ in the interior of $\Theta$ such that $\hat{\theta}^{ML} \overset{p}{\rightarrow} \theta^*$ and $\sqrt{n}(\hat{\theta}^{ML}-\theta^*)$ is asymptotically normal.
    \item[\textbf{\textbf{(A1)}}] The prior density $\pi(\theta)$ is continuous and positive in a neighborhood of $\theta^*$. 
    \item[\textbf{\textbf{(A2)}}] 
    Given $\Delta_{n,\theta^*}=\sqrt{n}(\hat{\theta}^{ML}-\theta^*)$, there exists a positive definite matrix $V_{\theta^*}$ such that
    \begin{equation}    
    \begin{split}    
    R_n(h) &\equiv \log\left(\frac{f_n(X^n|\theta^* + h/\sqrt{n})}{f_n(X^n|\theta^*)} \right) \\
    &\qquad - h^TV_{\theta^*}\Delta_{n,\theta^*} + \frac{1}{2}h^TV_{\theta^*}h,
    \end{split}
    \end{equation}
    satisfies $\sup_{h\in K}|R_n(h)|\rightarrow 0$ in $f_{0,n}$-probability for all compact sets $K \subseteq \mathbb{R}^p$.
    \item[\textbf{\textbf{(A3)}}] There exists a $\gamma > 0$ such that the following is finite
    \begin{equation*}    
    \mathbb{E}_{f_{0,n}} \left[\int_{\mathbb{R}^p} ||h||_2^{k_{0}(1+\gamma)}n^{-\frac{p}{2}}\pi_{n,\alpha} \left(\theta^* + \frac{h}{\sqrt{n}}\right)  dh\right],
    \end{equation*}
    for some $k_{0} \in \mathbb{N}$.
\end{itemize}
We adapt the proof of \cite[Theorem 1]{MORV22} -- which in turn adapts the arguments
 of \cite[Theorem 2.1]{KV12} -- to prove Theorem \ref{thm:moments_alt}. We retain  assumptions \textbf{(A0)} - \textbf{(A2)} from \cite{KV12} and \cite{MORV22}, but make a slightly stronger posterior concentration assumption in  \textbf{(A3)}. Indeed, our assumption implies the posterior concentration assumption in \cite{KV12} and \cite{MORV22}, which states that $\pi_{n,\alpha_n}(\cdot|X^n)$ concentrates at rate $\sqrt{n}$ around $\theta^*$ if for every sequence $r_n \rightarrow \infty$,
\begin{equation}\label{post_conc}  
\mathbb{E}_{f_{0,n}}\left[ \mathbb{P}_{\pi_{n,\alpha}(\theta|X^n)}\left(||\sqrt{n}(\theta - \theta^*)||_2 > r_n \right)\right] \rightarrow 0.
\end{equation}
Applying Markov's inequality to the left side of \eqref{post_conc} and using the change of variable $h = \sqrt{n}\left(\theta-\theta^*\right)$ gives us the upper bound
\begin{equation}\label{post_conc_markov}  
r_n^{-k_{0}}\mathbb{E}_{f_{0,n}}\left[ \int_{\mathbb{R}^p} ||h||_2^{k_{0}}n^{-\frac{p}{2}}\pi_{n,\alpha} \left(\theta^* + \frac{h}{\sqrt{n}}\right)  dh \right].
\end{equation}
Because $r_n \rightarrow \infty$, the posterior concentration assumption is satisfied if 
\begin{equation*}
\mathbb{E}_{f_{0,n}}\left[ \int_{\mathbb{R}^p} ||h||_2^{k_{0}}n^{-\frac{p}{2}}\pi_{n,\alpha} \left(\theta^* + \frac{h}{\sqrt{n}}\right) dh \right] < \infty,
\end{equation*}
which is implied by  \textbf{(A3)}. 

\subsection{Convergence of power posterior moments}\label{sec:moments}
We denote by  $\phi(\cdot|\mu,\Sigma)$ the density of a normal random vector with mean $\mu$ and covariance matrix $\Sigma$.
Recall that \cite{MORV22} showed the following result:
\begin{equation*} d_{TV}\left(\pi_{n,\alpha}\left(\theta|X^n\right),\phi\left(\theta \bigg| \hat{\theta}^{ML}, \frac{1}{\alpha n}V_{\theta^*}^{-1}\right)\right)\rightarrow 0,
\end{equation*}
in $f_{0,n}$-probability, where  $V_{\theta^*}$ is the matrix from assumption \textbf{(A2)} and it  can be thought of as measuring the curvature of the likelihood. Theorem \ref{thm:moments_alt} essentially shows that the $k^{th}$ moment of the $\alpha$-posterior approaches the $k^{th}$ moment of the limiting Gaussian from \cite{MORV22} with growing sample size.

\begin{Theorem}\label{thm:moments_alt}
Assume \textbf{\textbf{(A0)}}--\textbf{(A3)} hold. Then, for any integer $k\in[1,k_{0}]$,
\begin{equation*}
\begin{split}
\small
&    \int_{\mathbb{R}^p}\left|\left|h^{\otimes k}\right|\right|_1\left|\pi_{n,\alpha}^{LAN}(h) - \phi_n(h)\right|dh \rightarrow 0,
\end{split}
\end{equation*}
in $f_{0,n}$-probability, where
\begin{equation}\label{eq:tensor_norm}
    \left|\left|h^{\otimes k}\right|\right|_1 = \sum_{i_1=1}^{p} \cdots \sum_{i_k=1}^{p}\left|h_{i_1}\times\cdots\times h_{i_k}\right|,
\end{equation}
is the $1$-norm of a $k^{th}$ order tensor and
\begin{equation}
\begin{split}
&\pi_{n,\alpha}^{LAN}(h) \equiv n^{-p/2}\pi_{n,\alpha}(\theta^*+h/\sqrt{n}|X^n), \\
&\phi_n(h) \equiv n^{-p/2}\phi\left(h\big|\sqrt{n}(\hat{\theta}^{ML}-\theta^*),V_{\theta^*}^{-1}/\alpha\right),
\label{eq:scaled_densities}
\end{split}
\end{equation}
are the centered and scaled versions of the $\alpha$-posterior and its limiting normal. 
\end{Theorem}
\begin{Remark}
We mention that the objects in \ref{eq:scaled_densities} are proper densities. Indeed, using a change of variable $h = \sqrt{n}(\theta - \theta^*)$, we see
\begin{align*}
\int \pi_{n,\alpha}^{LAN}(h) dh &= \int n^{-p/2}\pi_{n,\alpha}(\theta^*+h/\sqrt{n}|X^n) dh \\
&= \int \pi_{n,\alpha}(\theta|X^n) d\theta = 1.
\end{align*}
Similarly,
\begin{equation*}
\begin{split}
\int\phi_n(h) dh &= \int n^{-\frac{p}{2}}\phi\left(h\big|\sqrt{n}(\hat{\theta}^{ML}-\theta^*), \frac{V_{\theta^*}^{-1}}{\alpha}\right) dh \\
&= \int \phi\left(\theta \big|\hat{\theta}^{ML}, \frac{V_{\theta^*}^{-1}}{n \alpha}\right)d\theta = 1.
\end{split}
\end{equation*}
\end{Remark}

\begin{Remark}
With a change of variable $h = \sqrt{n}(\hat{\theta}^{ML}-\theta^*)$, we see that as long as $k_0\geq 2$, Theorem~\ref{thm:moments_alt} shows that both of the following converge to zero in $f_{0,n}$-probability:
\begin{equation*}
\begin{split}
\small
&    \int_{\mathbb{R}^p} \hspace{-4pt} \sqrt{n}  \|\theta - \theta^*\|_1\left|\pi_{n,\alpha}(\theta| X^n) - \phi\bigg(\theta \bigg|\hat{\theta}^{ML}, \frac{V_{\theta^*}^{-1}}{n\alpha}\bigg)\right|d\theta,\\
&     \int_{\mathbb{R}^p} \hspace{-8.5pt}  n  \| \hspace{-1pt} [\hspace{-.5pt}  \theta \hspace{-2pt}  - \hspace{-2pt}  \theta^* \hspace{-.5pt}  ] \hspace{-.5pt}  [\hspace{-.5pt}  \theta \hspace{-2pt}  - \hspace{-2pt}  \theta^* \hspace{-.5pt}  ]^T \hspace{-1.5pt}  \|_1 \hspace{-2.5pt} \left|\hspace{-.5pt}\pi_{n,\alpha} \hspace{-.5pt} ( \hspace{-.5pt} \theta| X^n \hspace{-1pt}) \hspace{-1.5pt} - \hspace{-1.5pt}\phi \hspace{-1pt}  \bigg(\hspace{-2pt} \theta \bigg| \hat{\theta}^{ML} \hspace{-1pt}, \hspace{-1pt} \frac{V_{\theta^*}^{-1}}{n\alpha} \hspace{-1.5pt}   \bigg) \hspace{-1.5pt} \right|\hspace{-1.5pt} d\theta.
\end{split}
\end{equation*}
\end{Remark}

The proof of \cite{MORV22} and \cite{KV12} relates ``closeness" of the likelihood ratios of two densities to closeness in terms of total variation distance. Since the results in \cite{MORV22} and \cite{KV12} are Bernstein-von Mises results, we modify their proof to relate ``closeness" of likelihood ratios to ``closeness" of moments. To accomplish this, we modify Lemma 4 from \cite{MORV22} (Lemma \ref{lem:momentTVbound} in this work) to express the moments of the $\alpha$-posterior in terms of the moments and likelihood ratios of these two densities. Lemma \ref{lem:momentTVbound} reveals that the moments of the $\sqrt{n}$-scaled $\alpha$-posterior must be \textbf{(I)} finite on a growing neighborhood surrounding the pseudo-true parameter and Lemma \ref{MomentKc} shows that the moments are \textbf{(II)} vanishing on the complement of this neighborhood. Lemmas \ref{lem:momentTVbound} and \ref{MomentKc} require \textbf{(A3)} instead of the weaker posterior concentration assumption from \cite{MORV22}. 

\subsection{Asymptotic normality of the Bayes estimator}\label{sec:BayesEstimator}
The result in Section \eqref{sec:moments} is useful for determining the large sample behavior of point estimators derived from the $\alpha$-posterior. In particular, because $\alpha$ does not affect the centering of the $\alpha$-posterior, as the the mean of the posterior approaches $\hat{\theta}^{ML}$, it stands to reason that the limiting distribution of the mean of the $\alpha$-posterior is unaffected by $\alpha$. Our second result verifies this.

\begin{Theorem} \label{thm:BayesEstimator}
Assume \textbf{\textbf{(A0)}}--\textbf{\textbf{(A3)}} hold. Then the $\alpha$-posterior mean estimator $\hat{\theta}^B = \int_{\mathbb{R}^p}\theta\pi_{n,\alpha}(\theta|X^n)d\theta$ has the limiting distribution \small $\sqrt{n}(\hat{\theta}^B - \theta^*) \overset{d}{\rightarrow} N(0,\tilde{V}_{\theta^*}),$
\normalsize
where $\tilde{V}_{\theta^*}$ is the asymptotic variance of the MLE.
\end{Theorem}

\begin{Remark}
    If $\log f_n(X^n|\theta^*)$ is differentiable at $\theta^*$, then $\tilde{V}_{\theta^*} = V_{\theta^*}^{-1}M_{\theta^*}V_{\theta^*}^{-1}$, where $M_{\theta^*}$ is the limit in $f_{0,n}$-probability as $n\rightarrow\infty$ of $\frac{1}{n}\dot{\ell}_{n,\theta^*}\dot{\ell}_{n,\theta^*}^T$ and $\dot{\ell}_{n,\theta^*} = \nabla_{\theta}\log f_n(X^n|\theta)|_{\theta=\theta^*}$ is the score function of the model $f_n(X^n|\theta)$ at $\theta^*$ \cite{KV12}.
\end{Remark}

The high-level argument in our proof of Theorem \ref{thm:BayesEstimator} involves showing that the Bayes estimator is asymptotically equivalent to the (pseudo) MLE, which implies that they have the same limiting distribution. In other words, we need to show that $\sqrt{n}(\hat{\theta}^B-\hat{\theta}^{ML}) = o_p(1)$. Since $\hat{\theta}^B$ is defined as the first moment of the $\alpha$-posterior, we can obtain the desired result by examining the moments of the $\alpha$-posterior centered around the pseudo-true parameter, which we do in Theorem \ref{thm:moments_alt}.

\section{Discussion}
In this work, we have studied the moments of $\alpha$-posteriors in parametric, low-dimensional models to better understand point estimators derived from posterior tempering. We found that because posterior tempering does not affect the ``first-order" properties of the posterior, the limiting distribution of the posterior mean is also unaffected by tempering. This result demonstrates that while tempering may induce robustness in the posterior distribution, it does not do so in the mean of the tempered posterior. We believe this phenomenon should hold for other point estimators derived from posteriors (e.g., the posterior median and mode), which is left for future work.

We also mention some other directions for future work. First, as noted in Section \ref{sec:assumptions}, our \textbf{(A3)} is stronger than the posterior concentration assumption used to prove the Bernstein-von Mises results in \cite{MORV22, KV12}. However, we believe that posterior concentration is a sufficient condition for our Theorem \ref{thm:moments_alt} to hold for at least $k=1$ and $k=2$. Future work could involve modifying our proof to reflect this. Second, while we have advanced our understanding of $\alpha$-posteriors for fixed, constant, $\alpha$, we would like to establish similar guarantees to $\alpha$-posteriors where $\alpha$ is selected in a data-driven way (i.e., in settings where $\alpha$ might change with the sample size, $n$). 

\section{Proofs of Main Results}
\subsection{Proof of Theorem \ref{thm:moments_alt}}
The proof is adapted from the proofs of \cite[Theorem 1]{MORV22} and \cite[Theorem 2.1]{KV12} to show convergence of the posterior moments rather than the total variation distance. 

Let $Z_0$ denote the integral we are trying to control, namely
\begin{equation}
\label{eq:Z0_def}
Z_0 \coloneqq \int_{\mathbb{R}^p}\left|\left|h^{\otimes k}\right|\right|_1\left|\pi^{LAN}_{n,\alpha}(h) - \phi_n(h)\right|dh.
\end{equation}
Recalling the definition of $\left|\left|h^{\otimes k}\right|\right|_1$ in \ref{eq:tensor_norm} and using the fact that $\left|\left|h\right|\right|_1 \leq \sqrt{p}||h||_2$, we note that 
\begin{equation}\label{eq:tensor_norm_bound}
\begin{split}
    \left|\left|h^{\otimes k}\right|\right|_1 
    &= \sum_{i_1=1}^{p} \cdots \sum_{i_k=1}^{p}\left|h_{i_1}\times\cdots\times h_{i_k}\right| \\
    &= \sum_{i_1=1}^{p}\left|h_{i_1}\right|  \cdots \sum_{i_k=1}^{p}\left|h_{i_k}\right|= \left|\left|h\right|\right|_1^k \leq p^{\frac{k}{2}} \left|\left|h\right|\right|_2^k,
\end{split}
\end{equation}
with $p$ fixed. By the bound in \ref{eq:tensor_norm_bound} and monotonicity of the integral, we have that
\begin{equation}
\label{eq:Z_def}
Z_0 \leq p^{\frac{k}{2}}\int_{\mathbb{R}^p}\left|\left|h\right|\right|_2^k\left|\pi^{LAN}_{n,\alpha}(h) - \phi_n(h)\right|dh \coloneqq Z.
\end{equation}
Thus, it suffices to show that $\mathbb{E}_{f_{0,n}}[Z] \rightarrow 0$, which implies $\mathbb{E}_{f_{0,n}}[Z_0] \rightarrow 0$. This, in turn, implies that $Z_0 \overset{p}{\rightarrow} 0$ by Markov's inequality. We show convergence in $L_1$ by \textbf{(I)} showing convergence in $L_1$ conditioned on an intersection of events and \textbf{(II)} showing the complement of these events occur with vanishing probability. \\
\indent We now show convergence in $L_1$. We will use the following strategy to do this. For vectors $g,h\in K_0$, define the random variables
\begin{align}
& f_n^+(g,h)=\left\{1-\frac{\phi_n(h)\pi_{n,\alpha}^{LAN}(g|X^n) }{\pi_{n,\alpha}^{LAN}(h|X^n) \phi_n(g)}\right\}^+,  \label{fn+} \\
&f_n^-(g,h)=\left\{\frac{\pi_{n,\alpha}^{LAN}(h|X^n) \phi_n(g)}{\phi_n(h) \pi_{n,\alpha}^{LAN}(g|X^n)}-1\right\}^-,  \label{fn-}
\end{align}
where $\{x\}^+ = \max\{0,x\}$ and $\{x\}^-=\max\{0,-x\}$.
We appeal to \cite[Lemma 5]{MORV22} to show that $f_n^+(g,h)$ and $f_n^-(g,h)$ are well-defined with high-probability. Consider $\eta > 0$ and define the following events,
\begin{equation}
\begin{split}
\label{eq:events}
\mathcal{A} &= \left\{ \sup_{g,h\in B_{0}(r_n)}f_n^+(g,h) \leq \eta \right\}, \\
\mathcal{B} &= \left\{ \sup_{g,h\in B_{0}(r_n)}f_n^-(g,h) \leq \eta \right\}.
\end{split}
\end{equation}
The events $\mathcal{A}$ and $\mathcal{B}$ hold with high probability for large enough $n$ by \cite[Lemma 5]{MORV22}. To apply \cite[Lemma 5]{MORV22}, we mention that \textbf{(I)} because $\theta^*$ is in the interior of $\Theta \subseteq \mathbb{R}^p$ by \textbf{\textbf{(A1)}}, there exists a $\delta > 0$ such that ball $B_{\theta^*}(\delta) \subset \Theta$ and  $\pi$ is continuous and positive on the ball and \textbf{(II)} $\sqrt{n}$-stochastic local asymptotic normality holds by \textbf{(A2)}.\\
\indent Next, we upper bound $\mathbb{E}_{f_{0,n}}[Z]$ as follows using H{\"o}lder's Inequlity for some $\gamma > 0$:
\begin{equation}\label{Zbound}
\begin{split}
   & \mathbb{E}_{f_{0,n}}\hspace{-1pt}[Z] 
    \hspace{-1pt} = \hspace{-1pt} \mathbb{E}_{f_{0,n}} \hspace{-1pt} [Z\mathbbm{1}\{\mathcal{A} \cap \mathcal{B}\}]  + \hspace{-1pt} \mathbb{E}_{f_{0,n}} \hspace{-1pt} [Z\mathbbm{1}\{(\mathcal{A} \cap \mathcal{B})^c\}] \\
    &\leq \mathbb{E}_{f_{0,n}}[Z\mathbbm{1}\{\mathcal{A} \cap \mathcal{B}\}] \\
    &\quad + \mathbb{E}_{f_{0,n}}[Z^{1+\gamma}]^{\frac{1}{1+\gamma}}\left(\mathbb{P}_{f_{0,n}}(\mathcal{A}^c) + \mathbb{P}_{f_{0,n}}(\mathcal{B}^c)\right)^{\frac{\gamma}{1+\gamma}}.
\end{split}
\end{equation}
We aim to show that each term in \eqref{Zbound} converges to zero.\\
\noindent \textbf{First term in \eqref{Zbound}.}
We upper bound the first term of \eqref{Zbound}, which is $\mathbb{E}_{f_{0,n}}[Z\mathbbm{1}\{\mathcal{A} \cap \mathcal{B}\}]$.
Using \eqref{eq:Z_def} and appealing to Lemma~\ref{lem:momentTVbound} with $s(h) = p^{\frac{k}{2}}||h||_2^k$, we have that 
\small
\begin{equation}\label{lem3bd}
\begin{split}    
&Z = p^{\frac{k}{2}}\int_{\mathbb{R}^p}||h||_2^k\left|\pi_{n,\alpha}^{LAN}(h|X^n) - \phi_n(h)\right|dh \\
&\leq p^{\frac{k}{2}} \hspace{-2pt}\left[ \sup_{g,h\in B_{0}(r_n)}f^+_n(g,h)\right]\int_{\mathbb{R}^p}||h||_2^k\pi_{n,\alpha}^{LAN}(h|X^n)dh \\
&\quad + p^{\frac{k}{2}}\left[\sup_{g,h\in B_{0}(r_n)}f_n^-(g,h)\right]\int_{\mathbb{R}^p}||h||_2^k\phi_n(h)dh  \\
&\quad+ p^{\frac{k}{2}}\int_{||h||_2 > r_n}||h||_2^k \left(\pi_{n,\alpha}^{LAN}(h|X^n) +\phi_n(h)\right)dh.
\end{split}
\end{equation}
\normalsize
Therefore, to study $\mathbb{E}_{f_{0,n}}[Z\mathbbm{1}\{\mathcal{A} \cap \mathcal{B}\}]$, we notice that  the suprema in \eqref{lem3bd} are bounded. This gives the upper bound 
\begin{equation}\label{EZbd}
\begin{split}    
&\mathbb{E}_{f_{0,n}}[Z\mathbbm{1}\{\mathcal{A} \cap \mathcal{B}\}] \\
&\leq p^{\frac{k}{2}}\eta\mathbb{E}_{f_{0,n}}\left[\int_{\mathbb{R}^p}||h||_2^k\pi_{n,\alpha}^{LAN}(h|X^n)dh\right] \\
&\quad + p^{\frac{k}{2}}\eta\mathbb{E}_{f_{0,n}}\left[\int_{\mathbb{R}^p}||h||_2^k\phi_n(h)dh\right] \\
& \quad +p^{\frac{k}{2}}\mathbb{E}_{f_{0,n}}\left[\int_{||h||_2 > r_n}||h||_2^k \phi_n(h)dh \right]\\
&\quad+ p^{\frac{k}{2}}\mathbb{E}_{f_{0,n}}\left[\int_{||h||_2 > r_n}||h||_2^k \pi_{n,\alpha}^{LAN}(h|X^n) dh\right].
\end{split}
\end{equation}
Label the four terms in \eqref{EZbd} as $T_1$, $T_2$, $T_3$, and $T_4$. We argue that, for arbitrary $\epsilon > 0$, each can be  upper bounded by $\epsilon/4$ (possibly for large enough $n$).\\ 
\indent Consider term $T_1$ of \eqref{EZbd}. By \textbf{\textbf{(A3)}}, there exists $M_1 < \infty$ such that
\begin{equation*}    
\mathbb{E}_{f_{0,n}}\left[\int_{\mathbb{R}^p}||h||_2^k \, \pi_{n,\alpha}^{LAN}(h|X^n)dh\right] < M_1.
\end{equation*}
Thus, we can pick $\eta < \frac{\epsilon}{4M_1p^{k/2}}$ such that
\begin{equation*}    
T_1 = p^{k/2}\eta\mathbb{E}_{f_{0,n}}\left[\int_{\mathbb{R}^p}||h||_2^k  \, \pi_{n,\alpha}^{LAN}(h|X^n)dh\right] < \frac{\epsilon}{4}.
\end{equation*}
For term $T_2$, we notice that by a change of variables,
\begin{equation}
\begin{split}
\label{eq:t2_eq1}
&\int_{\mathbb{R}^p} \|h\|_2^k \phi_n(h) dh \\
&= \int_{\mathbb{R}^p} n^{-\frac{p}{2}}  \|h\|_2^k \phi\left(h\bigg|\sqrt{n}(\hat{\theta}^{ML}-\theta^*), \frac{V_{\theta^*}^{-1}}{\alpha}\right) dh \\
&= \int_{\mathbb{R}^p} n^{\frac{k}{2}} \|\theta - \hat{\theta}^{ML} + \theta^*\|_2^k\phi\left(\theta \bigg|0, \frac{V_{\theta^*}^{-1}}{n \alpha}\right)d\theta \\
&= n^{\frac{k}{2}} \mathbb{E}_{\mathcal{N}(0, V_{\theta^*}^{-1}/(n\alpha))} \left[ \|\theta - \hat{\theta}^{ML} + \theta^*\|_2^k \right].
\end{split}
\end{equation}
Next we upper bound the above using H{\"o}lder's Inequality\footnote{Note that for $m>0$ and  $v=(v_1,\dots,v_p)\in\mathbb{R}^p$, H{\"o}lder's inequality implies  $(\sum_{j=1}^p|v_j|)^m\leq p^{m-1}\sum_{j=1}^p|v_j|^m$} and centered absolute Gaussian moments:
\begin{equation}
\begin{split}
\label{eq:t2_eq2a}
&n^{\frac{k}{2}} \mathbb{E} \hspace{-2pt} \left[ \hspace{-1pt} \|\theta - \hat{\theta}^{ML} + \theta^*\|_2^k \right] \\
&= n^{\frac{k}{2}} \mathbb{E}\hspace{-2pt} \left[ \hspace{-2pt}\left( \sum_{i=1}^p (\theta_i - \hat{\theta}_i^{ML} + \theta_i^*)^2\hspace{-1pt}   \right)^{\frac{k}{2}}  \right]\\
&\leq n^{\frac{k}{2}} p^{\frac{k}{2} - 1} \sum_{i=1}^p \mathbb{E} \left[ (\theta_i - \hat{\theta}_i^{ML} + \theta_i^*)^k \right]\\
&\leq  n^{\frac{k}{2}} p^{\frac{k}{2} - 1} \sum_{i=1}^p \mathbb{E} \left[ \lvert\theta_i - \hat{\theta}_i^{ML} + \theta_i^*\lvert^k \right]
\\
&\leq n^{\frac{k}{2}} 2^{k-1} p^{\frac{k}{2} - 1} \sum_{i=1}^p \left(  \mathbb{E} \left[ \lvert\theta_i\lvert^k \right] +   \lvert \hat{\theta}_i^{ML} - \theta_i^*\lvert^k \right).
\end{split}
\end{equation}
The first term on the right side of \eqref{eq:t2_eq2a} equals
\begin{equation}
\begin{split}
\label{eq:t2_eq2b}
&n^{\frac{k}{2}} 2^{k-1} p^{\frac{k}{2} - 1} \sum_{i=1}^p \mathbb{E}[ \lvert\theta_i\lvert^k ] \\
&\quad  =  \frac{ 2^{\frac{3k}{2} - 1} p^{\frac{k}{2} - 1}}{\sqrt{\pi}}  \Gamma \left(\frac{k+1}{2}\right)  \sum_{i=1}^p  \left[\frac{[V_{\theta^*}^{-1}]_{ii}}{\alpha}\right]^{\frac{k}{2}},
\end{split}
\end{equation}
and the second term on the right side of \eqref{eq:t2_eq2a} equals $ 2^{k-1} p^{\frac{k}{2} - 1} \sum_{i=1}^p (\sqrt{n} \lvert \hat{\theta}_i^{ML} - \theta_i^*\lvert)^k$. Because $\sqrt{n}( \hat{\theta}_i^{ML} - \theta_i^*)$ is asymptotically normal by Assumption \textbf{(A0)}, there exists an $N_1 \coloneqq N_1(p, k, M)$ such that for all $n > N_1$, we have that the second term on the right side of \eqref{eq:t2_eq2a} is bounded by some $M < \infty$ with high probability. From \eqref{eq:t2_eq1} - \eqref{eq:t2_eq2b} we have shown for all $n > N_1$,
\begin{equation}
\begin{split} \label{eq:t2_eq3}
&\mathbb{E}_{f_{0,n}}\left[\int_{\mathbb{R}^p} \|h\|_2^k \phi_n(h) dh \right] \leq   M_2,
\end{split}
\end{equation}
where
\begin{equation*}
\begin{split} 
&M_2 \coloneqq  \frac{ 2^{\frac{3k}{2} - 1} p^{\frac{k}{2} - 1}}{\sqrt{\pi}}  \Gamma \left(\frac{k+1}{2}\right)  \sum_{i=1}^p  \left[\frac{[V_{\theta^*}^{-1}]_{ii}}{\alpha}\right]^{\frac{k}{2}}+ M.
\end{split}
\end{equation*}
Putting this together, we have
\begin{equation*}
\begin{split}
&T_2 = p^{\frac{k}{2}}\eta \mathbb{E}_{f_{0,n}}\left[\int_{\mathbb{R}^p} \|h\|_2^k \phi_n(h) dh \right] \leq p^{\frac{k}{2}}\eta M_2.
\end{split}
\end{equation*}
Hence, by choosing $\eta < {\epsilon }/{(4 M_2 p^{\frac{k}{2}})}$ we have $T_2 \leq \epsilon/4$.

We control $T_3$ by leveraging \eqref{eq:t2_eq3}. Indeed, by Cauchy-Schwarz,
\begin{equation}\label{eq:T3bd2}
\begin{split}
    &T_3 
    = p^{\frac{k}{2}}\mathbb{E}_{f_{0,n}}\left[ \int_{\mathbb{R}^p} \|h\|_2^k(1-\mathbbm{1}\{||h||_2\leq r_n\}) \phi_n(h) dh\right] \\
    &\leq p^{\frac{k}{2}}\mathbb{E}_{f_{0,n}}\Bigg[ \sqrt{\int_{\mathbb{R}^p} \|h\|_2^{2k} \phi_n(h) dh} \\
    &\qquad\quad\times \sqrt{\int_{\mathbb{R}^p} (1-\mathbbm{1}\{||h||_2\leq r_n\}) \phi_n(h) dh} \, \Bigg] \\
    &\leq p^{\frac{k}{2}}\sqrt{\mathbb{E}_{f_{0,n}}\Bigg[ \int_{\mathbb{R}^p} \|h\|_2^{2k} \phi_n(h) dh\Bigg]} \\
    &\qquad\quad\times \sqrt{\mathbb{E}_{f_{0,n}}\Bigg[\int_{\mathbb{R}^p} (1-\mathbbm{1}\{||h||_2\leq r_n\}) \phi_n(h) dh\Bigg]}.
\end{split}
\end{equation}
By appealing to \eqref{eq:t2_eq3} with the power $2k$ instead of $k$, we know that $\mathbb{E}_{f_{0,n}}[ \int_{\mathbb{R}^p} \|h\|_2^{2k} \phi_n(h) dh] \leq  \tilde{M}_2$ for $n > N_1$.
Because $r_n\rightarrow\infty$ and $\phi_n(h) \leq (2\pi)^{-p/2}|V_{\theta^*}^{-1}/\alpha|^{-1/2}$ for all $h\in\mathbb{R}^p$, we see that $(1-\mathbbm{1}\{||h||_2\leq r_n\})\phi_n(h) \rightarrow 0$. Thus, by the bounded convergence theorem, there exists $N_2\coloneqq N_2(\epsilon, \tilde{M}_2, p, k)$ such that for all $n > N_2$,
\[
\mathbb{E}_{f_{0,n}}\Bigg[\int_{\mathbb{R}^p} (1-\mathbbm{1}\{||h||_2\leq r_n\}) \phi_n(h) dh\Bigg] < \frac{\epsilon^2}{16\tilde{M}_2p^{k}}.
\]
Thus, for $n > \max(N_1, N_2)$, we have that $T_3 < \frac{\epsilon}{4}$.\\
\indent By appealing to Lemma \ref{MomentKc} (assumptions are met by appealing to \textbf{\textbf{(A3)}}) with $f_Z(z)$ being $\pi_{n,\alpha}^{LAN}(z|X^n)$, we pick $N_3\coloneqq  N_3(p, k, \epsilon, \gamma)$ such that for all $n > N_3$,
\begin{equation*}    
\mathbb{E}_{f_{0,n}} \hspace{-3pt} \left[\int_{||h||_2 > r_n} \hspace{-7pt} ||h||_2^k\pi_{n,\alpha}^{LAN}(h|X^n)dh\right] \hspace{-2pt} < \hspace{-1pt} \frac{\epsilon}{4p^{k/2}}.
\end{equation*}
Thus,
\begin{equation*}    
T_4 = p^{\frac{k}{2}}\mathbb{E}_{f_{0,n}} \hspace{-3pt} \left[\int_{||h||_2 > r_n} \hspace{-7pt} ||h||_2^k\pi_{n,\alpha}^{LAN}(h|X^n)dh\right] \hspace{-2pt} < \hspace{-1pt} \frac{\epsilon}{4}
\end{equation*}
for $n > N_3$. Therefore, we have shown that the first term of \eqref{Zbound} is upperbounded by $\epsilon$ whenever $\eta < \frac{\epsilon}{4 \min\{ M_1, M_2\} p^{\frac{k}{2}}}$ 
and $n \geq \max\{N_1, N_2, N_3\}$.

\noindent \textbf{Second term in \eqref{Zbound}.}
First, we consider upper bounding the term $\mathbb{E}_{f_{0,n}}[Z^{1+\gamma}]$. 
Notice that by the fact that $|\pi_{n,\alpha}^{LAN}(h|X^n) - \phi_n(h)| \leq \pi_{n,\alpha}^{LAN}(h|X^n) + \phi_n(h)$,
\begin{equation*}
\begin{split}
    &Z^{1+\gamma} = \left[p^{\frac{k}{2}}\hspace{-2pt}\int_{\mathbb{R}^p}\hspace{-1pt}||h||_2^k\left|\pi_{n,\alpha}^{LAN}(h|X^n) - \phi_n(h)\right|dh\right]^{1+\gamma} \\
    &\leq \Bigg[2p^{\frac{k}{2}}\max\bigg\{\int_{\mathbb{R}^p}   ||h||_2^k\pi_{n,\alpha}^{LAN}(h|X^n)dh, \\
    &\hspace{100pt} \int_{\mathbb{R}^p}  ||h||_2^k \phi_n(h)dh\bigg\}\Bigg]^{1+\gamma},
\end{split}
\end{equation*}
Next, applying Jensen's Inequality, we have that
\begin{equation*}
\begin{split}
    &\Bigg[\hspace{-3pt}\max \hspace{-2pt} \bigg\{ \hspace{-3pt} \int_{\mathbb{R}^p}  \hspace{-4pt} ||h||_2^k\pi_{n,\alpha}^{LAN}\hspace{-2pt}(h|X^n)dh,  \int_{\mathbb{R}^p} \hspace{-4pt}  ||h||_2^k \phi_n(h)dh \hspace{-2pt} \bigg\} \hspace{-2pt} \Bigg]^{1+\gamma} \\
    &= \max\Bigg\{\bigg(\int_{\mathbb{R}^p}   ||h||_2^k\pi_{n,\alpha}^{LAN}(h|X^n)dh\bigg)^{1+\gamma}, \\
    &\hspace{100pt} \bigg(\int_{\mathbb{R}^p}   ||h||_2^k \phi_n(h)dh\bigg)^{1+\gamma}\Bigg\} \\
    &\leq  \max\Bigg\{ \int_{\mathbb{R}^p} ||h||_2^{k(1+\gamma)}\pi_{n,\alpha}^{LAN}(h|X^n)dh,\\
    &\hspace{100pt} \int_{\mathbb{R}^p}  ||h||_2^{k(1+\gamma)} \phi_n(h)dh\Bigg\} \\
    &\leq \int_{\mathbb{R}^p}  ||h||_2^{k(1+\gamma)} (\pi_{n,\alpha}^{LAN}(h|X^n) + \phi_n(h)) dh.
\end{split}
\end{equation*}
Therefore, we have the bound
\begin{equation*}
\begin{split}
    & \mathbb{E}_{f_{0,n}}\hspace{-1pt}[Z^{1+\gamma}] \leq (2p^{\frac{k}{2}})^{1+\gamma} \mathbb{E}_{f_{0,n}}\hspace{-2pt} \left[\int_{\mathbb{R}^p} \hspace{-2pt} ||h||_2^{k(1+\gamma)} \phi_n(h)dh\right]\\
     &\hspace{20pt} + \hspace{-2pt} (2p^{\frac{k}{2}})^{1+\gamma}\mathbb{E}_{f_{0,n}} \hspace{-2.5pt} \left[\int_{\mathbb{R}^p} \hspace{-2pt}||h||_2^{k(1+\gamma)}\pi_{n,\alpha}^{LAN}\hspace{-1pt} (h|X^n)dh  \right]\hspace{-2pt}.
\end{split}
\end{equation*}
We argue that the two terms on the right side of the above are upper bounded by constants.
By Assumption \textbf{\textbf{(A3)}}, there exists an $M_3 < \infty$ such that $\mathbb{E}_{f_{0,n}}[\int_{\mathbb{R}^p}||h||_2^{k(1+\gamma)}\pi_{n,\alpha}^{LAN}(h|X^n)dh] < M_3$. By appealing to \eqref{eq:t2_eq3} with $k(1 + \gamma)$ instead of $k$, we know that for $n > N_1$,
\begin{equation*}
\mathbb{E}_{f_{0,n}}\left[\int_{\mathbb{R}^p}||h||_2^{k(1+\gamma)}\phi_n(h)dh\right] < M_4.
\end{equation*}
Thus, we have shown that
\begin{equation*}
\begin{split}
    \mathbb{E}_{f_{0,n}}[Z^{1+\gamma}] 
    &\leq (2p^{k/2})^{1+\gamma}(M_3 + M_4) \\
    &\leq 2(2p^{k/2})^{1+\gamma} \max\{M_3, M_4\}.
\end{split}
\end{equation*}

We finally upper bound $\mathbb{P}(\mathcal{A}^c)$ and $\mathbb{P}(\mathcal{B}^c)$.
First, note that $f_n^+(g,h)$ is the same as $f_n(g,h)$ from \cite{MORV22}. Next, we use the fact that $\left\{1-x\right\}^+=\left\{x-1\right\}^-$ to see
\begin{equation}\label{eq:fn_equality}
\begin{split}
    f_n^-(g,h)
    &=\left\{\frac{\pi_{n,\alpha}^{LAN}(h|X^n) \phi_n(g)}{\phi_n(h) \pi_{n,\alpha}^{LAN}(g|X^n)}-1\right\}^-\\
    &=\left\{1-\frac{\pi_{n,\alpha}^{LAN}(h|X^n) \phi_n(g)}{\phi_n(h) \pi_{n,\alpha}^{LAN}(g|X^n)}\right\}^+
    =f_n^+(h,g).
\end{split}
\end{equation}
Taking the supremum in $g,h$ over $B_0(r_n)$ of both sides of \ref{eq:fn_equality} gives us that
\begin{equation}\label{eq:fn_sup_equality}
    \sup_{g,h\in B_{0}(r_n)}f_n^-(g,h)=\sup_{g,h\in B_{0}(r_n)}f_n^+(h,g).
\end{equation}
Thus, there exists $N_4\coloneqq N_4(\eta, \epsilon, M_3, M_4, \gamma, p, k)$ such that $\mathbb{P}(\mathcal{A}^c) = \mathbb{P}(\mathcal{B}^c)$ is sufficiently small for all $n>N_4$ by \cite[Lemma 5]{MORV22}:
\begin{equation*}
\begin{split}
\mathbb{P}(\mathcal{A}^c)=\mathbb{P}(\mathcal{B}^c)
&< \frac{\epsilon^{\frac{1+\gamma}{\gamma}}}{2^{1+\frac{2+\gamma}{\gamma}} p^{\left(\frac{k}{2}\right)\left(\frac{1+\gamma}{\gamma}\right)}\hspace{-1.5pt}\max\{M_3, M_4\}^{\frac{1}{\gamma}}}.
\end{split}
\end{equation*}
Putting this together, we find that for $n> N_4$,
\begin{align*}
&\mathbb{E}_{f_{0,n}}[Z^{1+\gamma}]^{\frac{1}{1+\gamma}}\left(\mathbb{P}_{f_{0,n}}(\mathcal{A}^c) + \mathbb{P}_{f_{0,n}}(\mathcal{B}^c)\right)^{\frac{\gamma}{1+\gamma}}\\
&\leq  [2(2p^{k/2})^{1+\gamma}  \max\{M_3, M_4 \}]^{\frac{1}{1+\gamma}}  \\
&\quad\times \left[ \frac{\epsilon^{\frac{1+\gamma}{\gamma}}}{2^{\frac{2+\gamma}{\gamma}} p^{\left(\frac{k}{2}\right)\left(\frac{1+\gamma}{\gamma}\right)}\max\{M_3, M_4\}^{\frac{1}{\gamma}}}  \right]^{\frac{\gamma}{1+\gamma}} = \epsilon.
\end{align*}

\noindent \textbf{Final bound.}
We have shown that whenever $\eta$ is chosen small enough, for all $n > \max(N_1, N_2, N_3, N_4)$,  for arbitrary $\epsilon > 0$,
\begin{equation*}
\begin{split}    
\mathbb{E}[Z] 
&= p^{k/2}\mathbb{E}\Bigg[ \int_{\mathbb{R}^p}  ||h||_2^{k}\left|\pi_{n,\alpha}^{LAN}(h|X^n) - \phi_n(h)\right|dh\Bigg] \\
&< 2\epsilon.
\end{split}
\end{equation*}
This gives us the desired result. 

\subsection{Proof of Theorem \ref{thm:BayesEstimator}}
Recall, $\hat{\theta}^B = \int_{\mathbb{R}^p}\theta\pi_{n,\alpha}(\theta|X^n)d\theta$. Using the change of variable $h = \sqrt{n}(\theta-\theta^*)$,
    \begin{equation*}
    \begin{split}    
    &\hat{\theta}^B = \int\theta\pi_{n,\alpha}(\theta|X^n)d\theta = \frac{1}{\sqrt{n}}\int h\pi_{n,\alpha}^{LAN}(h)dh + \theta^*.
    \end{split}
    \end{equation*}
    In the above, we used the definition $\pi_{n,\alpha}^{LAN}(h)= n^{-\frac{p}{2}} \pi_{n,\alpha}(\theta^* + \frac{h}{\sqrt{n}}|X^n)$ in \eqref{eq:scaled_densities}. Recalling the definition of $\phi_n(h)$ in \eqref{eq:scaled_densities}, we can write,
    \begin{equation}\label{normalized_bayes}
    \begin{split}    
        &\sqrt{n}(\hat{\theta}^B   - \theta^*) = \int h\pi_{n,\alpha}^{LAN}(h)dh  \\
        &= \int h(\pi_{n,\alpha}^{LAN}(h) - \phi_n(h))dh + \int h\phi_n(h)dh \\
        &= \int h(\pi_{n,\alpha}^{LAN}(h) - \phi_n(h))dh  + \sqrt{n}(\hat{\theta}^{ML} \hspace{-2pt}  - \theta^*).
    \end{split}
    \end{equation}
    The second term of \ref{normalized_bayes} converges weakly to  $\mathcal{N}(0,\tilde{V}_{\theta^*})$ by Assumption \textbf{(A0)}. It remains to show that the first term of the last line of \ref{normalized_bayes} is $o_p(1)$. This follows from 
    \begin{equation*}
    \begin{split}    
        &\left|\left| \int h(\pi_{n,\alpha}^{LAN}(h) - \phi_n(h))dh \right|\right|_1 \\
        &\leq \int \left|\left|h\right|\right|_1\left|\pi_{n,\alpha}^{LAN}(h) - \phi_n(h)\right|dh = o_p(1)
    \end{split}
    \end{equation*}
    where the equality follows from Theorem \ref{thm:moments_alt} for $k=1$.
\section{Technical Lemmas}
\begin{Lemma}\label{lem:momentTVbound}
Consider sequences of densities $\varphi_n$ and $\psi_n$. For a given compact set $K \subset \mathbb{R}^p$, suppose that the densities $\varphi_n$ and $\psi_n$ are positive on $K$. Then for any function, $s(h)$, that is nonnegative on all of $\mathbb{R}^p$,
\begin{align*}
    &\int_{\mathbb{R}^p}s(h)\left|\varphi_n(h) - \psi_n(h)\right|dh \\
    &\leq \left[ \sup_{g,h\in K}\tilde{f}^+_n(g,h)\right]\int_{\mathbb{R}^p}s(h)\psi_n(h)dh \\
    &\quad + \left[\sup_{g,h\in K}\tilde{f}_n^-(g,h)\right]\int_{\mathbb{R}^p}s(h)\varphi_n(h)dh \\
    &\quad + \int_{\mathbb{R}^p\setminus K}s(h) \left(\psi_n(h) + \varphi_n(h) \right)dh,
\end{align*}
where
\begin{equation*}
\begin{split}
\tilde{f}_n^+(g,h) \coloneqq \left\{ 1 - \frac{\varphi_n(h)\psi_n(g)}{\psi_n(h)\varphi_n(g)}\right\}^+, \\
\tilde{f}_n^-(g,h) \coloneqq \left\{\frac{\psi_n(h)\varphi_n(g)}{\varphi_n(h)\psi_n(g)}-1\right\}^-.
\end{split}
\end{equation*}
\end{Lemma}
\begin{Proof}
    Since $s(h)$, $\psi_n(h)$, and $\varphi_n(h)$ are nonnegative on all of $\mathbb{R}^p$, we have that $s(h)|\psi_n(h) - \varphi_n(h)| \leq s(h)(\psi_n(h) + \varphi_n(h))$. Using this fact and that $|\psi_n(h) - \varphi_n(h)| = \{\psi_n(h) - \varphi_n(h)\}^- + \{\psi_n(h) - \varphi_n(h)\}^+$,
    \begin{equation}\label{rhs}
    \begin{split}    
    &\int_{\mathbb{R}^p}s(h)|\psi_n(h) - \varphi_n(h)| dh \\
    &= \int_{K}s(h)|\psi_n(h) - \varphi_n(h)|dh \\
    &\quad + \int_{\mathbb{R}^p\setminus K}s(h)|\psi_n(h) - \varphi_n(h)|dh \\
    & \leq \int_{K}s(h)\{\psi_n(h) - \varphi_n(h)\}^+dh \\
    & \quad +  \int_{K}s(h)\{\psi_n(h) - \varphi_n(h)\}^-dh \\
    & \quad + \int_{\mathbb{R}^p\setminus K}s(h)(\psi_n(h) + \varphi_n(h))dh.
    \end{split}
    \end{equation}
     We bound the first two terms on the right side of \eqref{rhs} separately. Define $a_n = (\int_K \psi_n(g)dg)^{-1}$ and $b_n = (\int_K\varphi_n(g)dg)^{-1}$ and assume $a_n \geq b_n$ without loss of generality. For the first term in \eqref{rhs}, for all $h\in K$,
    \[
    \frac{\varphi_n(h)}{\psi_n(h)} = \frac{a_n}{b_n}\int_K\frac{\varphi_n(h)}{\psi_n(h)}\frac{\psi_n(g)}{\varphi_n(g)}b_n\varphi_n(g)dg.
    \]
    Thus, 
    \begin{equation}\label{rhs1}
    \begin{split}    
    &\left\{\psi_n(h)-\varphi_n(h)\right\}^+ = \left\{1 - \frac{\varphi_n(h)}{\psi_n(h)}\right\}^+\psi_n(h) \\
    &= \left\{1-\frac{a_n}{b_n}\int_K\frac{\varphi_n(h)}{\psi_n(h)}\frac{\psi_n(g)}{\varphi_n(g)}b_n\varphi_n(g)dg\right\}^+\psi_n(h). 
    \end{split}
    \end{equation}
    Applying Jensen's inequality on the convex function $f(x) =\{ 1-x\}^+$, we have $\{ 1 - \mathbb{E}[X]\}^+ \leq \mathbb{E}[\{ 1 - X\}^+]$. Applying this to \eqref{rhs1},
    \begin{equation} \label{rhs2}
    \begin{split}
    &\left\{1-\frac{a_n}{b_n}\int_K\frac{\varphi_n(h)}{\psi_n(h)}\frac{\psi_n(g)}{\varphi_n(g)}b_n\varphi_n(g)dg\right\}^+\\
    &\leq \int_K\left\{1-\frac{a_n}{b_n}\frac{\varphi_n(h)}{\psi_n(h)}\frac{\psi_n(g)}{\varphi_n(g)}\right\}^+b_n\varphi_n(g)dg \\
    &= \int_K \tilde{f}_n^+(g,h)b_n\varphi_n(g)dg.
    \end{split}
    \end{equation}
    The final step above uses that when $a_n/b_n \geq 1$, we have $\{ 1-(a_n/b_n)x\}^+ \leq \{1-x\}^+$ for $x \geq 0$. 
    Using \eqref{rhs1}-\eqref{rhs2} we can rewrite the first term of \eqref{rhs} as
    \begin{align*}
    &\int_{K}s(h)\{\psi_n(h) - \varphi_n(h)\}^+dh \\
    & \leq \int_K\int_K s(h) \tilde{f}^+_n(g,h)b_n\varphi_n(g)\psi_n(h)dgdh \\
    &\leq \hspace{-2pt} \left[ \sup_{g,h\in K}\tilde{f}^+_n(g,h)\right] \hspace{-2pt} \int_K \hspace{-2pt} s(h)\psi_n(h)\left[\int_K b_n\varphi_n(g)dg\right] \hspace{-2pt}dh\\
    &\leq \hspace{-2pt} \left[ \sup_{g,h\in K}\tilde{f}^+_n(g,h)\right] \hspace{-2pt} \int_K \hspace{-2pt} s(h)\psi_n(h) \hspace{-2pt}dh,\\
    &\leq \hspace{-2pt} \left[ \sup_{g,h\in K}\tilde{f}^+_n(g,h)\right] \hspace{-2pt} \int_{\mathbb{R}^p} \hspace{-2pt} s(h)\psi_n(h) \hspace{-2pt}dh.
    \end{align*}
    We turn to the second term of \eqref{rhs} and write 
    \begin{equation*}
    \begin{split}
    &\big\{\psi_n(h)-\varphi_n(h)\big\}^- = \left\{\frac{\psi_n(h)}{\varphi_n(h)}-1\right\}^-\varphi_n(h)\\
    &= \left\{\frac{a_n}{b_n}  \int_K  \frac{\psi_n(h) \varphi_n(g)}{\varphi_n(h) \psi_n(g)}b_n\psi_n(g)dg  - 1\right\}^-  \varphi_n(h). 
    \end{split}
    \end{equation*}
    Applying Jensen's inequality to the convex function $f(x)=\{x-1\}^-$, we have that $\{\mathbb{E}[X]-1\}^- \leq \mathbb{E}[\{X-1\}^-]$ and
    \begin{equation*}
    \begin{split}    
    & \left\{\frac{a_n}{b_n}\int_K\frac{\psi_n(h) \varphi_n(g)}{\varphi_n(h) \psi_n(g)}b_n\psi_n(g)dg-1  \right\}^-  \\
    &\leq \int_K \hspace{-3pt} \left\{ \frac{a_n\psi_n(h) \varphi_n(g) }{b_n \varphi_n(h) \psi_n(g)} - 1  \right\}^- b_n\psi_n(g)dg  \\
    &\overset{(a)}\leq  \int_K  \left\{\frac{\psi_n(h)\varphi_n(g)}{\varphi_n(h) \psi_n(g)}-1\right\}^- b_n\psi_n(g)dg \\
    &=  \int_K \tilde{f}_n^-(g,h)b_n\psi_n(g)dg \overset{(b)}\leq  \int_K \tilde{f}_n^-(g,h)a_n\psi_n(g)dg.
     \end{split}
    \end{equation*}
     Step $(a)$ uses that $\{ (a_n/b_n)x-1\}^-\leq \{x-1\}^-$ for $x \geq 0$ when $a_n/b_n \geq 1$  and step $(b)$ that $b_n \leq a_n$. Hence, we can rewrite  the second term of \eqref{rhs} as 
      \begin{equation*}
    \begin{split}  
    & \int_{K}s(h)\{\psi_n(h) - \varphi_n(h)\}^-dh \\
    & \leq \int_K\int_K  s(h) \tilde{f}_n^-(g,h)a_n\psi_n(g)\varphi_n(h)dgdh \\
    &\leq \left[\sup_{g,h\in K}\hspace{-2pt}\tilde{f}_n^-(g,h)\right]\hspace{-2pt}\int_K s(h)\varphi_n(h) \hspace{-2pt}\left[\int_K \hspace{-2pt} a_n\psi_n(g)dg\right] \hspace{-2pt}dh.
    \end{split}
    \end{equation*}
    This establishes the the final bound since
    \begin{equation*}
    \begin{split}
    &\int_K s(h) \varphi_n(h) \left[\int_K  a_n\psi_n(g)dg\right]dh \\
    &= \int_K s(h) \varphi_n(h)dh \leq \int_{\mathbb{R}^p}s(h)\varphi_n(h)dh.
    \end{split}
    \end{equation*}
\end{Proof}

\begin{Lemma}\label{MomentKc}
    Suppose a random variable $Z \overset{d}{\coloneqq} Y|X^n$ has a density $f_Z(\cdot)$ on $\mathbb{R}^p$. Assume there exists a $\gamma > 0$ such that $\mathbb{E}_{f_{0,n}}[\mathbb{E}[||Z||_2^{k(1+\gamma)}]] < \infty$. Then, for any $\epsilon > 0$ and $r_n \rightarrow \infty$, there exists an integer $N(\epsilon, \gamma, k) > 0$ such that for all $n > N(\epsilon, \gamma, k)$, 
    \begin{align}
    &\mathbb{E}_{f_{0,n}}\left[\int_{||z||_2 > r_n}||z||_2^k \, f_Z(z) dz\right] < \epsilon, \label{eq:lemma3_eq1}
    \end{align}
\end{Lemma}
\begin{Proof}
For an arbitrary sequence $r_n\rightarrow\infty$, we choose $N_0$ large enough that $r_n > 0$ for all $n > N_0$. Let $z = (z_1,\ldots,z_p)\in\mathbb{R}^p$ and let $f_{||z||_2^k}(\cdot)$ denote the density function of $||z||_2^k$. Define the random variable 
\begin{equation*}
\begin{split}
    V 
    &\coloneqq ||Z||_2^k\mathbbm{1}\left\{ ||Z||_2 > r_n\right\}
    = \begin{cases}
    0 &\textrm{if } ||Z||_2 \leq r_n, \\
    ||Z||_2^k &\textrm{if } ||Z||_2 > r_n.
    \end{cases}
\end{split}
\end{equation*}
Thus, the density of $V$ is given by
\begin{equation*}
    f_V(v) =
    \begin{cases}
    \mathbb{P}(||Z||_2 \leq r_n) &\textrm{ if } v = 0, \\
    f_{||z||_2^k}(v) &\textrm{ if } v > r_n.
    \end{cases}
\end{equation*}
Next, notice that
\begin{equation}\label{markov0}
    \begin{split}
        \mathbb{E}[V]
        &= \int_0^{r_n} \mathbb{P}(V > t)dt + \int_{r_n}^\infty \mathbb{P}(V > t)dt.
    \end{split}
\end{equation}

Using the definition of $V$ and its density function, we study \eqref{markov0}.
Consider the first term on the right side,
\begin{equation}\label{markov}
    \begin{split}
        \int_0^{r_n}   \hspace{-1pt} \mathbb{P}(V > t)dt  &= \int_0^{r_n} \int_t^\infty f_V(v)dvdt \\
        &\overset{(a)}{=} \int_0^{r_n} \mathbb{P}(||Z||_2^k > r_n)dt  \\
        &= r_n\mathbb{P}(||Z||_2^k > r_n) \overset{(b)}{\leq} \frac{ \mathbb{E}[||Z||_2^{k(1+\gamma)}]}{r_n^\gamma} .
    \end{split}
\end{equation}
Step $(a)$ in \eqref{markov} follows because, for $0 < t < r_n$,
\begin{equation*}
\int_t^\infty f_V(v)dv = \int_{r_n}^\infty f_{||z||_2^k}(v) dv = \mathbb{P}(||Z||_2^k > r_n),
\end{equation*}
and step $(b)$ by Markov's inequality. 

For the second term on the right side of \eqref{markov0}, we again use Markov's Inequality as
\begin{equation}\label{markov1}
    \begin{split}
       \int_{r_n}^\infty \mathbb{P}(V > t)dt &\leq   \int_{r_n}^\infty \frac{\mathbb{E}[V^{1+\gamma}]}{t^{1+\gamma}} dt \\
       &\overset{(c)}{\leq} \int_{r_n}^\infty    \frac{\mathbb{E}[||Z||_2^{k(1+\gamma)}]}{t^{1+\gamma}} dt   \\
       &=   \frac{\mathbb{E}[||Z||_2^{k(1+\gamma)}]}{\gamma r_n^{\gamma}}.
       \end{split}
\end{equation}
Step $(c)$ in \eqref{markov1} follows because
\begin{equation*}
\begin{split}
    &V^{1 + \gamma}  = \left(||Z||_2^k\mathbbm{1}\left\{||Z||_2 > r_n\right\}\right)^{1 + \gamma}  \leq ||Z||_2^{k(1+\gamma)}.    
\end{split}
\end{equation*}

Plugging \eqref{markov} and \eqref{markov1} into \eqref{markov0}, we find
\begin{equation}\label{markov_new}
    \begin{split}
        \mathbb{E}[V]
        &\leq \frac{\gamma + 1}{\gamma r_n^\gamma}  \mathbb{E}\left[||Z||_2^{k(1+\gamma)}\right] .
    \end{split}
\end{equation}
Taking the expectation with respect to $X^n$ in \eqref{markov_new} gives,
\begin{equation}\label{eq:final_bound}
\begin{split}
    \mathbb{E}_{f_{0,n}}\left[\mathbb{E}[V]\right] &= \mathbb{E}_{f_{0,n}}\left[\mathbb{E}[||Z||_2^k\mathbbm{1}\left\{ ||Z||_2 > r_n\right\}]\right] \\
    &\leq  \frac{\gamma + 1}{\gamma r_n^\gamma}  \mathbb{E}_{f_{0,n}}\left[\mathbb{E}[||Z||_2^{k(1+\gamma)}]\right] < \epsilon,
\end{split}
\end{equation}
where the final bound  in the equations above follows since $\mathbb{E}_{f_{0,n}}\left[\mathbb{E}[||Z||_2^{k(1+\gamma)}]\right] < \infty$ by assumption and $r_n\rightarrow\infty$. Thus, there exists $N(\epsilon, \gamma, k) \geq N_0$ such that the bound is true for all $n > N(\epsilon, \gamma, k)$, which proves \eqref{eq:lemma3_eq1}.
\end{Proof}
\section*{Acknowledgment}
This research was partially supported by the NSF grant DMS-231097 (M.\ Avella Medina).
\printbibliography
\end{document}